\documentclass[12pt]{article}

\UseRawInputEncoding

\usepackage{mathtools,amsthm,amssymb,color,cases}
\usepackage{mathrsfs}
\usepackage[abbrev]{amsrefs}
\usepackage{tikz}
\usepackage{enumitem}
\usepackage{setspace}

\usepackage{latexsym}
\usepackage{mathrsfs}
\usepackage{comment}
\usepackage{microtype}

\usepackage{indentfirst}

\usepackage{geometry}
\usepackage{etoolbox}
\AtEndEnvironment{proof}{\setcounter{proofpart}{0}}

 \newtheoremstyle{mystyle}
    {}
    {}
    {\normalfont}
    {}
    {\bfseries}
    {}
    { }
    {}
  \theoremstyle{mystyle}     

\allowdisplaybreaks[1]
\theoremstyle{definition}
\newtheorem{thm}{Theorem}[section]

\newtheorem{rem}[thm]{Remark}
\newtheorem*{rem*}{Remark}

\newtheorem{definition}{Definition}

\newtheorem{prop}[thm]{Proposition}

\newtheoremstyle{part}{}{}{\normalfont}{}{\itshape}{.}{.5em}{}
\theoremstyle{part}

\newcommand\blfootnote[1]{%
  \begingroup
  \renewcommand\thefootnote{}\footnote{#1}%
  \addtocounter{footnote}{-1}%
  \endgroup
}

\numberwithin{equation}{section}
\numberwithin{thm}{section}

\newenvironment{equ*}{
    \begin{equation*}
}{
    \end{equation*}
}

\newenvironment{pf}
   {{\noindent \bf Proof.}}{\hfill \qed}

\newtheoremstyle{part}{}{}{\normalfont}{}{\itshape}{.}{.5em}{}
\theoremstyle{part}

\newcommand{\ddj}{\dot{\Delta}_j}

\newcommand{\e}{\epsilon}

\newcommand{\pt}{\partial}

\newcommand{\F}{\mathcal {F}}

\newcommand{\Z}{\mathbb{Z}}

\newcommand{\R}{\mathbb R}

\newcommand{\intr}{\int_{\R^3}}

\newcommand{\al}{\alpha}

\newcommand{\x}{\xi}

\renewcommand{\div}{{\rm {div}}}

\newcommand{\supp}{\text{ supp }}

\renewcommand{\P}{\mathcal{P}}
\newcommand{\Q}{\mathcal{Q}}

\newcommand*{\dd}{\mathop{}\!\mathrm{d}}

\newcommand{\ha}{\hat{a}}
\newcommand{\hv}{\hat{v}}
\newcommand{\hp}{\hat{\phi}}
\newcommand{\hps}{\hat{\Psi}}

\DeclarePairedDelimiter{\norm}{\lVert}{\rVert}
\DeclarePairedDelimiter{\bignorm}{\bigg\|}{\bigg\|}

\DeclarePairedDelimiter{\bigc}{\big{(}}{\big{)}}
\DeclarePairedDelimiter{\Bigc}{\Big{(}}{\Big{)}}
\DeclarePairedDelimiter{\Biggc}{\Bigg{(}}{\Bigg{)}}

\DeclarePairedDelimiter{\Bigf}{\Big{[}}{\Big{]}}



\begin{document}

\begin{center}
    {\bf \large Optimality of the Decay Estimate of Solutions to the Linearised
Curl-Free Compressible Navier-Stokes Equations} \\
    
    \phantom{}
    
    Tsukasa Iwabuchi* \quad D\'aith\'i \'O hAodha** \\
    
     \phantom{}
    
    Mathematical Institute, Tohoku University,
    
    980-0845
\end{center}

\blfootnote{Email: *t-iwabuchi@tohoku.ac.jp, **david.declan.hughes.p6@dc.tohoku.ac.jp}

\begin{center}
\begin{minipage}{135mm}
\footnotesize
{\sc Abstract. }
We discuss optimal estimates of solutions to the compressible Navier-Stokes equations in Besov norms. In particular, we consider the estimate of the curl-free part of the solution to the linearised equations, in the homogeneous case.
We prove that our estimate is optimal in the $L^\infty$-norm by showing that the norm is bounded from below by the same decay rate. 
\end{minipage}
\end{center}

\section{Introduction}

This paper is concerned with the barotropic compressible Navier-Stokes system
\begin{align} \label{CNSo}
    \begin{cases}
    \pt_t \rho + \div (\rho u) = 0, & \text{ in } (0,\infty) \times\R^3, \\
    \pt_t (\rho u) + \div (\rho u \otimes u) - \div(2 \mu D(u) + \lambda \div (u) \text{Id}) + \nabla p = 0, & \text{ in } (0,\infty) \times\R^3, \\
    (\rho, u) |_{t=0} = (\rho_0, u_0), & \text{ in } \R^3,
    \end{cases}
\end{align}
where $\rho : [0,\infty) \times \R^3 \to [0,\infty),$ and $u: [0,\infty) \times \R^3 \to \R^3$ are unknown functions, representing the density and velocity of a fluid, respectively. 
$p: [0,\infty) \times \R^3 \to \R$ is the pressure in the fluid, and the barotropic assumption gives us $p \coloneqq p(\rho)$, for some smooth function $p(\cdot)$.
$\mu,\lambda$ are viscosity coefficients, taken such that
\begin{align*}
    \mu > 0, \quad 2\mu + \lambda > 0.
\end{align*}
We define the deformation tensor
\begin{align*}
    D(u) \coloneqq \frac{1}{2} \Big{(} Du + Du^T \Big{)}.
\end{align*}


In this paper, we examine the large-time behaviour, particularly time-decay estimates, of solutions to a linearised version of the above problem.
Matsumura-Nishida showed in~\cites{matsumura-nishida1979, matsumura-nishida1980} that \eqref{CNSo} has global solutions when equipped with data $(\rho_0, u_0)$ that is a small perturbation in $L^1 \cap H^3$ of $(\bar{\rho}, 0)$ for any positive constant $\bar{\rho}$, and proved the following decay result
\begin{align*}
    \bignorm{ 
    \begin{bmatrix}
        \rho(t) - \bar{\rho} \\
        u(t)
    \end{bmatrix}  
    }_{2}
    \leq C ({1+t})^{-3/4}.
\end{align*}
This is the decay rate of the solution to the heat equation with initial data in $L^1$.
Ponce then extended these results to other $L^p$ norms. In particular, for $p \in [2,\infty]$, $k\in\{0,1,2\}$, and dimension $d=2,3$,
\begin{align*}
    \bignorm{ 
    \nabla^k \begin{bmatrix}
        \rho(t) - \bar{\rho} \\
        u(t)
    \end{bmatrix}  
    }_{p}
    \leq C ({1+t})^{-\frac{d}{2}(1-\frac{1}{p}) - \frac{k}{2} }.
\end{align*}

Our results contribute to the theory of Compressible Navier-Stokes in Besov spaces. Global existence of strong solutions to \eqref{CNSo} for initial data $(\rho_0, u_0)$ in critical Besov spaces $ \dot{B}^{d/2}_{2,1} \times \dot{B}^{d/2 - 1}_{2,1} $ was first proven by Danchin in~\cite{Danchin2000} and large-time estimates in Besov norms for $p$ close to $2$ were proven by Danchin-Xu in~\cite{Danchin-Xu2017}.

Our goal in this paper is to obtain an optimal bound of the solution to a linearised version of the system \eqref{CNSo} which isolates the curl-free part of the velocity.
In particular, we prove a bound from above, in terms of time $t$, of the norm of solutions over space $x$.
We then prove that said bound is optimal (sharp) via a corresponding bound from below. 
Our proof is entirely self-contained.
The bound from below in our main result, Theorem \ref{main result}, is original, and is the first ever proof that the decay rate in $L^\infty$ is optimal, or sharp.
The estimates from above are analogous to those found in \cites{hoff-zumbrun, kobayashi-shibata}. 
Indeed, in \cite{hoff-zumbrun}, the following estimate is obtained.
\begin{prop} (\cite{hoff-zumbrun})
    Let $m\coloneqq \rho u$, $m_0 \coloneqq \rho_0 u_0$. Assume that 
    \[
    E \coloneqq 
    \bignorm{  
    \begin{bmatrix}
        a_0 \\
        m_0
    \end{bmatrix}
    }_1
    +
    \bignorm{
    \begin{bmatrix}
        a_0 \\
        m_0
    \end{bmatrix}}_{H^{ 1 + l }}
    \]
    is sufficiently small, where $l\geq3$ is an integer.
Then the Navier-Stokes system \eqref{CNSo} with initial data $\rho_0, u_0$ has a global solution satisfying the following decay estimate for all  for any multi-index $\al$ with $|\al| \leq (l-3)/2$:
\begin{align}
        \bignorm{
        D^\al_x \Bigc{
        \begin{bmatrix}
         a(t) \\
         m(t) - e^{t\Delta} \P m_0
        \end{bmatrix}}
        }_p
        \leq
        C(l) E (1+t)^{ - \frac{3}{2} (1-\frac{1}{p})  - \frac{1}{2} ( 1 - \frac{2}{p} )   - \frac{|\al|}{2}   }.
\end{align}
\end{prop}
Note that, in the norm in the above inequality, we are removing the divergence-free part of the linear term of $m$. Thus what remain are the nonlinear term and the curl-free part of the linear term.

Kobayashi-Shibata in \cite{kobayashi-shibata} also obtain the above decay rate. In fact, they obtained a refined estimate that requires less regularity on the initial data and separates the solution into high and low frequencies (see Definition \ref{besov spaces} below).
The above decay rate found in \cite{hoff-zumbrun} is associated with the low-frequency part of solutions, while the high-frequency part decays exponentially with $t$.

In what follows, we shall assume that the density approaches $1$ at infinity; and so we are concerned with strong solutions which are small perturbations from a constant state $(\rho,u) \equiv (1,0)$. 
We shall also assume that $\mu, \lambda$ are constant,
and set $a \coloneqq \rho - 1$. 
Our system \eqref{CNSo} can thus be rewritten into the following linearised problem:
\begin{align*}
    \begin{cases}
    \pt_t a + \div(u) = f & \text{ in } (0,\infty) \times \R^3, \\
    \pt_t u - \mu \Delta u - (\lambda + \mu) \nabla \div(u) + P'(1)\nabla a = g & \text{ in } (0,\infty) \times \R^3, \\
    (a,u)\Big{|}_{t=0} = (a_0,u_0) & \text{ in }\R^3,
    \end{cases}
\end{align*}
for some functions $f,g.$
We apply the orthogonal projections $\mathcal{P}$ and $\mathcal{Q}$ onto the divergence and curl-free fields, respectively. Then, setting $\al \coloneqq P'(1)$ and $\nu \coloneqq \lambda + 2\mu$, we get the system
\begin{align} \label{div and curl free split}
    \begin{cases}
    \pt_t a + \div (\mathcal{Q} u) = f & \text{ in } \R_{> 0} \times \R^3, \\
    \pt_t \Q u - \nu \Delta \Q u + \al \nabla a = \Q g & \text{ in } \R_{> 0} \times \R^3, \\
    \pt_t \P u - \mu \Delta \P u = \P g & \text{ in } \R_{> 0} \times \R^3.
    \end{cases}
\end{align}
We set 
\[
v \coloneqq |D|^{-1} \div ( u), \text{ where }
|D|^s u  \coloneqq \F^{-1} \Big{[} |\xi|^s \hat{u}\Big{]}, \ s\in\R.
\]
We note that one can obtain $v$ from $\Q u$ by a Fourier multiplier of homogeneous degree zero. Thus, bounding $v$ is equivalent to bounding $\Q u$ (see \cite{Danchinbook}, Lemma 2.2).

We note that we can set $\al = \nu = 1,$ without loss of generality, since the following rescaling
\begin{align*}
    a(t,x) = \tilde{a}\Bigc{\frac{\al}{\nu} t, \frac{\sqrt{\al}}{\nu} x }, \quad u(t,x) = \sqrt{\al} \ \tilde{u} \Bigc{ \frac{\al}{\nu} t , \frac{\sqrt{\al}}{\nu} x }
\end{align*}
ensures that $(\tilde{a}, \tilde{v})$ solves \eqref{div and curl free split} with $\al = \nu = 1.$
Thus we get that $(a, v)$ solves the following system:
\begin{align} \label{lin system av}
    \begin{cases}
    \pt_t a + |D| v = f & \text{ in } \R_{> 0} \times \R^3, \\
    \pt_t v - \Delta v - |D| a = h \coloneqq |D|^{-1} \div (g) & \text{ in } \R_{> 0} \times \R^3.
    \end{cases}
\end{align}

In this paper, we shall focus on the homogeneous case, where $f = h = 0$.
Before we give our main result, we introduce the concepts of Besov spaces, and high-frequency and low-frequency norms.
We denote by $P(\R^3)$ the set of all polynomials over $\R^3$.

\begin{definition} \label{besov spaces}
We use the Littlewood-Paley decomposition of unity to define homogeneous Besov spaces. Let $\{ \hat{\phi}_j \}_{ j \in\Z}$ be a set of non-negative measurable functions such that 
\begin{enumerate}
    \item $\displaystyle \sum_{ j \in\Z} \hat{\phi}_j (\x) = 1, \text{ for all } \x \in \R^3 \backslash \{0\}$,
    \item $\hat{\phi}_j (\x) = \hat{\phi}_0(2^{-j}\x)$,
    \item $\supp \hat{\phi}_j (\x) \subseteq \{ \x \in \R^3 \ | \ 2^{j-1} \leq |\x| \leq 2^{j+1} \}$.
\end{enumerate}
For $f \in \mathcal{S'}(\R^3)/P(\R^3)$, we write 
\[
\dot{\Delta}_j f \coloneqq \F^{-1} [\hat{\phi}_j \hat{f}].
\]
The Besov norm is then defined as follows: for $1 \leq p,q \leq \infty$, and $s \in \R$, we define 
\[
\norm{f}_{\dot{B}^{s}_{p,q}} \coloneqq \Bigc{ \sum_{j \in \Z} 2^{sqj} \norm{\dot{\Delta}_j f}^q_{p} }^{1/q}.
\]
The set $\dot{B}^{s}_{p,q}$ is defined as the set of functions, $f \in \mathcal{S'}(\R^3)/P(\R^3)$, whose Besov norm is finite.
Finally, we also regularly use the following notation for so-called high and low-frequency norms:
\[
\norm{f}^h_{\dot{B}^{s}_{p,q}} \coloneqq \Bigc{ \sum_{j \geq 3} 2^{sqj} \norm{\dot{\Delta}_j f}^q_{p} }^{1/q}, \quad
\norm{f}^l_{\dot{B}^{s}_{p,q}} \coloneqq \Bigc{ \sum_{j \leq 2} 2^{sqj} \norm{\dot{\Delta}_j f}^q_{p} }^{1/q}.
\]
\end{definition}

We now give our main result. We obtain the same decay rate as the previous proposition, and also prove that this decay rate is optimal in the $L^\infty$-norm.
\begin{thm} \label{main result}
Let $s \in \R$, $p \in [2,\infty],$ $q \in [1,\infty]$, and $t > 1$.
Then
\begin{align} \label{thm1 l1 linf}
        \bignorm{ \F^{-1} 
        \begin{bmatrix}
        \ha(t) \\
        \hv(t)
        \end{bmatrix}  }_{\dot{B}^s_{p,q}}
        \leq 
        C t^{-\frac{3}{2} (1-\frac{1}{p}) - \frac{1}{2} (1 - \frac{2}{p}) }
        \bignorm{ \F^{-1} 
        \begin{bmatrix}
        \ha_0 \\
        \hv_0
        \end{bmatrix}  }^l_{\dot{B}^{s}_{1,q}}
        +
        C e^{-t}
        \bignorm{ \F^{-1} 
        \begin{bmatrix}
        \ha_0 \\
        \hv_0
        \end{bmatrix}  }^h_{
        \dot{B}^{s+3(1-1/p)}_{1,q} }.
\end{align}
For the high-frequency norm, we also have for all $t>0$,
\begin{align} \label{thm1 lp lp}
        \bignorm{ \F^{-1} 
        \begin{bmatrix}
        \ha(t) \\
        \hv(t)
        \end{bmatrix}  }^h_{\dot{B}^s_{p,q}}
        \leq 
        C e^{-t}
        \bignorm{ \F^{-1} 
        \begin{bmatrix}
        \ha_0 \\
        \hv_0
        \end{bmatrix}  }^h_{\dot{B}^{s}_{p,q} }.
\end{align}
Also, there exists initial data such that, for all sufficiently large $t$,
\begin{align} \label{thm 1 linf below}
    \bignorm{ \F^{-1} \begin{bmatrix}
        \ha(t) \\
        \hv(t)
        \end{bmatrix} }_{\infty} \geq C t^{-2}.
\end{align}
\end{thm}

\begin{rem} \label{remark 1}
    The derivative loss in the high-frequency part of the $L^1 \textendash L^\infty$ estimate \eqref{thm1 l1 linf} is in fact caused by only one of the two eigenvalues (derived in Section 3) associated with the fundamental solution of our problem. 
    The other eigenvalue has strong enough smoothing properties that no additional derivatives are required when estimating its associated semigroup (see Proposition \ref{mainprop}).
    Note, however, that we can avoid the need for more regular initial data by using \eqref{thm1 lp lp} instead for the high-frequency part of solutions.
\end{rem}

\bigskip

We offer some commentary on the above decay rate for the low frequency part.
In Hoff-Zumbrun~\cite{hoff-zumbrun}, it is demonstrated that the curl-free part of the solution to \eqref{CNSo} is approximated by a convolution of the Green's function for a parabolic system (whence we obtain the heat-like decay) and a function which decays like the fundamental solution of a wave equation (whence the additional decay).
To show that this decay rate is sharp in $L^\infty$, we find that this wave-like part of the decay is slowest along $|x| = t$.
This is a natural result. 
Indeed, looking at the homogeneous case $f=h=0$ for \eqref{lin system av}, we may write
\begin{align*}
    & \pt^2_t a - \Delta a - |D|^3 v =0,
    & \pt^2_t v - \Delta v - \pt_t \Delta v =0.
\end{align*}
Note that the wave propagation speed above is $1$. This is a consequence of our assumption made on the constant $\al$.

\section{Preliminaries}














\begin{definition}(The Fourier Transform)
For a function, $f$, we define the Fourier transform of $f$ as follows:
\begin{align*}
\F[f](\x) \coloneqq \hat{f}(\x) \coloneqq \frac{1}{(2\pi)^{3/2}}\intr e^{-i x\cdot\x} f(x) \dd x.
\end{align*}
The inverse Fourier transform is then defined as 
\begin{align*}
\F^{-1}[\hat{f}](x) \coloneqq \frac{1}{(2\pi)^{3/2}} \intr e^{i x \cdot \xi} \hat{f}(\x) \dd \x.
\end{align*}
For the purpose of calculating inequalities, we will frequently omit the factor of $1/(2\pi)^{3/2}$.
\end{definition}




\begin{definition}(Orthogonal Projections on the divergence and curl-free fields)
The projection mapping $\P$ is a matrix with each component defined as follows for $i,j \in \{1, 2, 3\}$:
\begin{align*}
(\P)_{i,j} \coloneqq \delta_{i,j} + (-\Delta)^{-1} \pt_i \pt_j.
\end{align*}
We then define $\Q \coloneqq 1 - \P$.
For $f \in (\dot{B}^s_{p,q}(\R^3))^3$, with $s \in \R$, and $p, q \in [1,\infty]$, we may write
\begin{align*}
\P f \coloneqq (1 + (-\Delta)^{-1} \nabla \div) f.
\end{align*}
\end{definition}

\section{Proof of Main Result}

We begin our analysis of $(a,v)$ which solves the system
\begin{align*}
    \begin{cases} 
    \pt_t a + |D| v = 0 & \text{ in } (0,\infty) \times \R^3, \\
    \pt_t v - \Delta v - |D| a = 0 & \text{ in } (0,\infty) \times \R^3.
    \end{cases}
\end{align*}
Taking the Fourier transform over space $x$, we can write the above system as
\begin{align} \label{ODE}
    \frac{\dd }{\dd t}
    \begin{bmatrix}
    \ha \\
    \hv
    \end{bmatrix}
    =
    M_{|\x|} \begin{bmatrix}
    \ha \\
    \hv
    \end{bmatrix}
    ,
    \quad \text{with}
    \quad M_{|\x|} \coloneqq 
    \begin{bmatrix}
    0 & -|\x| \\
    |\x| & -|\x|^2
    \end{bmatrix}.
\end{align}
Then we may write the following formula for the solution to \eqref{ODE}:
\begin{align*}
    \begin{bmatrix}
    \ha (t) \\
    \hv (t)
    \end{bmatrix} 
    = 
    e^{t M_{|\x|}} 
    \begin{bmatrix}
    \ha_0 \\
    \hv_0
    \end{bmatrix}
    .
\end{align*}

We obtain the following eigenvalues for $M_{|\x|}$, which differ between high and low frequencies:
\begin{align*}\lambda_{\pm}(\x) \coloneqq 
    \begin{cases}
    -\frac{|\x|^2}{2} \Bigc{1 \pm i { 
    \sqrt{\frac{4}{|\x|^2} - 1 } } }, & \text{ for } |\x| <2, \\
    -\frac{|\x|^2}{2} \Bigc{1 \pm { 
    \sqrt{ 1 - \frac{4}{|\x|^2} } } }, & \text{ for } |\x| >2,
    \end{cases}
\end{align*}
from which we obtain the following formulas for our solution:
\begin{align}
    \label{a formula}
    \ha(t) & = \frac{e^{t\lambda_-}\lambda_+ - e^{t\lambda_+ } \lambda_- }{ \lambda_+ - \lambda_- }\ha_0 
    + \frac{ \bigc{ e^{t\lambda_-} - e^{t\lambda_+} } \lambda_+ }{ \lambda_+ - \lambda_- }\hv_0,
    \\
    \label{v formula}
    \hv(t) & = \frac{ \bigc{ e^{t\lambda_+} - e^{t\lambda_-} } \lambda_- }{ \lambda_+ - \lambda_- }\ha_0 + 
    \frac{ e^{t\lambda_+} \lambda_+ - e^{t\lambda_-} \lambda_- }{ \lambda_+ - \lambda_- }\hv_0.
\end{align}
We note that, for high and low frequencies away from $|\x| = 2$, the eigenvalues can be bounded above by a mere constant.
However, both the numerator and the denominator approach $0$ as $|\x| \to 2$.
We thus must consider how to deal with the $L^2$ and $L^\infty$ norms of 
\[
\frac{e^{t\lambda_-}\lambda_+ - e^{t\lambda_+ } \lambda_- }{ \lambda_+ - \lambda_- } \hp_j
\]
and the other similar expressions in \eqref{a formula} and \eqref{v formula}, for $j \in \{ 0,1,2 \}$.
For details, see Appendix \ref{appendix no singularity}.
Otherwise, it suffices to prove time-decay estimates for the semigroup $e^{t\lambda_\pm}$ without including the eigenvalues in the $L^\infty$ norm for all $j\in\Z$, and in the $L^2$ norm for all $j \notin \{0,1,2\}$.

We will proceed with an interpolation argument after bounding the $L^\infty$ and $L^2$ norms explicitly. 
The $L^2$ norm estimate is identical to that of the heat semigroup, as no additional decay is obtained from the oscillation.
The $L^\infty$ norm yields a faster decay, for which we now give the proof.

\begin{prop} \label{mainprop} ($L^1 \textendash L^\infty$ estimates for high and low frequencies)
For all $t > 1$,
\begin{align} \label{inf one}
    \bignorm{\F^{-1} \Big{[} e^{ t \lambda_+}
    \hat{\phi}_j
    \Big{]} }_\infty
    & 
    \leq
    C e^{-t}, \text{ for all } j \geq 3,
    \\
    \bignorm{\F^{-1} \Big{[} e^{ t \lambda_-}
    \hat{\phi}_j
    \Big{]} }_\infty
    & \leq C 2^{3j}  e^{-t} , \text{ for all } j \geq 3, \label{inf two}
    \\
    \bignorm{\F^{-1} \Big{[} e^{ t \lambda_\pm}
    \hat{\phi}_j
    \Big{]} }_\infty
    & 
    \leq
    C t^{-2}, \text{ for all } j\leq 2. \label{inf three}
\end{align}
\end{prop}

\begin{rem}
    As we mentioned in Remark \ref{remark 1}, it is only the estimate of $e^{t\lambda_-}$ in \eqref{inf two} that causes a derivative loss.
\end{rem}

\begin{pf}
We start with \eqref{inf one}. We have $j \geq 3$, and so
\begin{align*}
    \supp \hp_j \subseteq \{ \x \in \R^3 \ | \ |\x| > 2 \},
\end{align*}
and thus we are dealing with the high-frequency definitions of our eigenvalues $\lambda_\pm$. That is, we have
\begin{align}
    \bignorm{\F^{-1} \Big{[} e^{ t \lambda_+} \notag
    \hp_j
    \Big{]} }_\infty
    & = 
    \bignorm{ \F^{-1} \Big{[} e^{-t\frac{|\x|^2}{2} \big{(} 1 + \sqrt{ 1 - \frac{4}{|\x|^2}} \big{)} }
    \hat{\phi}_j \Big{]} }_{\infty} 
    \\
    & = \notag
    \sup_{x \in \R^3} \Big{|}
    \intr
    e^{ix\cdot\x}
    e^{- t\frac{|\x|^2}{2} \sqrt{ 1 - \frac{4}{|\x|^2}} }
    e^{-t\frac{|\x|^2}{2}}
    \hat{\phi}(2^{-j} \x)
    \dd\x
    \Big{|}
    \\
    & \leq C \notag
    \intr
    \Big{|}
    e^{-t\frac{|\x|^2}{2}}
    \hp (2^{-j} \x)
    \Big{|}
    \dd \x
    \\
    & \leq C e^{-t} \notag
    \int_{ \{ \x \, | \,  |\x| > 2 \} }
    e^{-t|\x|^2/4}
    \dd \x
    \\
    & \leq C e^{-t}, \text{ for all } j \geq 3, \ t >1. \label{ high frequency og }
\end{align}

For \eqref{inf two}, the proof is slightly different. 
Before we continue, note that
    \begin{align*}
        {-t\frac{|\x|^2}{2} \bigc{ 1 - \sqrt{1 - {4}/{|\x|^2}}  }}
        =
        {-2t\bigc{1 + \sqrt{1 - {4}/{|\x|^2}}}^{-1}}
        =
        - t - \frac{4t}{|\x|^2} \bigc{1 +\sqrt{1 - 4 / |\x|^2}}^{-2}.
    \end{align*}
Thus, we may estimate as follows,
\begin{align*}
    \bignorm{ \F^{-1} \Big{[} e^{-t\frac{|\x|^2}{2} \big{(} 1 - \sqrt{ 1 - \frac{4}{|\x|^2}} \big{)} }
    \hat{\phi}_j \Big{]} }_{\infty}
    & \leq 
    \intr
    e^{-2t \bigc{ 1 + \sqrt{1 - \frac{4}{|\x|^2}} }^{-1}}
    \hat{\phi}_j (\x)
    \dd\x
    \\
    & \leq C 2^{3j} 
    e^{-t}
    .
\end{align*}

Now we consider \eqref{inf three}. For frequencies near $|\x|=2$, i.e. $0 \leq j \leq 2$, an exponential decay follows similarly to the above.
Finally, we are left with the low frequencies $j \leq -1$.
We start by extracting the heat-like decay by a change of variables.
\begin{align}
    & \sup_{x\in\R^3} \Big{|} 
    \intr
    e^{ix\cdot\x}
    e^{\pm it\frac{|\x|^2}{2} \sqrt{ \frac{4}{|\x|^2} - 1 }}
    e^{-t\frac{|\x|^2}{2}}
    \hp(2^{-j}\x)
    \dd\x
    \Big{|} \notag
    \\
    & = t^{-3/2}
    \sup_{x\in\R^3} \Big{|} 
    \intr
    e^{i (t^{-1/2}x)\cdot\x}
    e^{\pm i t^{1/2} \frac{|\x|^2}{2} \sqrt{ \frac{4}{|\x|^2} - t^{-1} }}
    e^{-\frac{|\x|^2}{2}}
    \hp(2^{-j} t^{-1/2} \x )
    \dd\x
    \Big{|}. \label{above integral}
\end{align}
Thus it remains to extract a further decay of $t^{-1/2}$ from the $L^\infty$ norm above.
We may estimate \eqref{above integral} 
to obtain the desired decay for all sufficiently large $t$:
\begin{align*}
    \sup_{x\in\R^3} \Big{|} 
    \intr
    e^{i (t^{-1/2}x)\cdot\x}
    e^{\pm i t^{1/2} \frac{|\x|^2}{2} \sqrt{ \frac{4}{|\x|^2} - t^{-1} }}
    e^{-\frac{|\x|^2}{2}}
    \hp(2^{-j} t^{-1/2} \x )
    \dd\x
    \Big{|} 
    & \leq C t^{ - 1/2}.
\end{align*}
We explicitly write the proof of the above in Appendix \ref{appendix prop 1}.
\end{pf}

\begin{prop} \label{proof of 1st}
    Let $s \in \R$, $p \in [2,\infty],$ and $q \in [1,\infty]$. 
    Then, for all $t>1$,
    \begin{align*}
        \bignorm{
        \F^{-1} \Big{[}
        e^{t \lambda_+}
        \hat{f}
        \Big{]}
        }_{\dot{B}^s_{p,q}}
        & \leq
        C t^{-\frac{3}{2} (1-\frac{1}{p}) - \frac{1}{2} (1 - \frac{2}{p}) }
        \norm{ f }^l_{\dot{B}^s_{1,q}} +
        C e^{-t}
        \norm{ f }^h_{\dot{B}^{s}_{1,q}},
        \\
        \bignorm{
        \F^{-1} \Big{[}
        e^{t \lambda_-}
        \hat{f}
        \Big{]}
        }_{\dot{B}^s_{p,q}}
        & \leq
        C t^{-\frac{3}{2} (1-\frac{1}{p}) - \frac{1}{2} (1 - \frac{2}{p}) }
        \norm{ f }^l_{\dot{B}^s_{1,q}} +
        C e^{-t} 
        \norm{ f }^h_{\dot{B}^{s+3(1-1/p)}_{1,q}}.
    \end{align*}
\end{prop}


\begin{pf}
    We start with the $L^\infty$ norm. Note that, by Proposition \ref{mainprop}, we have for all $j \geq 3$
    \begin{align*}
        \bignorm{
        \F^{-1} \Big{[}
        e^{t \lambda_+}
        \hp_j
        \hat{f}
        \Big{]}
        }_{\infty}
        & \leq
        C e^{-t}
        \norm{\dot{\Delta}_j f }_1,
        \\
        \bignorm{
        \F^{-1} \Big{[}
        e^{t \lambda_-}
        \hp_j
        \hat{f}
        \Big{]}
        }_{\infty}
        & \leq
        C 2^{3j}
        e^{-t}
        \norm{\dot{\Delta}_j f }_1.
    \end{align*}
    Next, for $j\leq 2$,
    \begin{align*}
        \bignorm{
        \F^{-1} \Big{[}
        e^{t \lambda_\pm}
        \hp_j
        \hat{f}
        \Big{]}
        }_{\infty}
        & \leq
        \bignorm{
        \F^{-1} \Big{[}
        e^{t \lambda_\pm}
        \hat{\Phi}_j
        \Big{]}
        }_{\infty}
        \norm{\dot{\Delta}_j f}_1
        \\
        & \leq
        C t^{-\frac{3}{2} - \frac{1}{2} }
        \norm{\dot{\Delta}_j f}_1,
    \end{align*}
    where $\hat{\Phi}_j \coloneqq \hp_{j-1} + \hp_j + \hp_{j+1}$. 
    
Next the $L^2$ norm for $j\leq2$ follows simply from Young's convolution inequality and heat kernel estimates:
    \begin{align*}
        \bignorm{
        \F^{-1} \Big{[}
        e^{t \lambda_\pm}
        \hp_j
        \hat{ f }
        \Big{]}
        }_{2}
        & \leq
        \bignorm{
        \F^{-1} \Big{[}
        e^{t \lambda_\pm}
        \hat{{\Phi}}_j
        \Big{]}
        }_{2}
        \norm{\dot{\Delta}_j f }_1
        \\
        & \leq
        C t^{-\frac{3}{4} }
        \norm{\dot{\Delta}_j f }_1.
    \end{align*}
For $j \geq 2,$ we get
\begin{align*}
    \bignorm{
        \F^{-1} \Big{[}
        e^{t \lambda_+}
        \hp_j
        \hat{ f }
        \Big{]}
        }_{2}
        & \leq
        \bignorm{
        \F^{-1} \Big{[}
        e^{t \lambda_+}
        \hat{\Phi}_j
        \Big{]}
        }_{2}
        \norm{\dot{\Delta}_j f }_1
        \\
        & \leq
        C e^{-t}
        \norm{\dot{\Delta}_j f }_1,
    \\
    \bignorm{
        \F^{-1} \Big{[}
        e^{t \lambda_-}
        \hp_j
        \hat{ f }
        \Big{]}
        }_{2}
        & \leq
        \bignorm{
        \F^{-1} \Big{[}
        e^{t \lambda_-}
        \hat{{\Phi}}_j
        \Big{]}
        }_{2}
        \norm{\dot{\Delta}_j f }_1
        \\
        & \leq
        C 2^{\frac{3}{2} j}  e^{-t}
        \norm{\dot{\Delta}_j  f }_1.
\end{align*}
    
    We use H\"older interpolation for $2<p<\infty$ to obtain $L^p$ norm estimates. Multiplying both sides by $2^s$ and taking the $\ell^q$ norm completes the proof.
\end{pf}


\begin{prop} \label{proof of 2nd}
($L^p \textendash L^p$ estimates for high frequencies) 
There exists a constant $C$ such that, for all $p \in [1,\infty]$ and $t>0$,
\begin{align*}
    \bignorm{
    \F^{-1} \Big{[}
    e^{t\lambda_\pm}
    \hp_j
    \hat{f}
    \Big{]}
    }_p
    & 
    \leq
    C e^{-t} \norm{  f }_p,  \text{ for all } 
    j \geq  3.
\end{align*}
\end{prop}

\noindent\textbf{Proof.} 
First, we make use of the same equalities as in Proposition \ref{mainprop}. That is, 
    \begin{align*}
        {-t\frac{|\x|^2}{2} \bigc{ 1 \pm \sqrt{1 - {4}/{|\x|^2}}  }}
        =
        {-2t\bigc{1 \mp \sqrt{1 - {4}/{|\x|^2}}}^{-1}}
        =
        - t - \frac{4t}{|\x|^2} \bigc{1 \mp \sqrt{1 - 4 / |\x|^2}}^{-2}.
    \end{align*}
Thus, we can rewrite
    \begin{align}
\notag    \bignorm{
    \F^{-1} \Big{[}
    e^{t\lambda_\pm}
    \hp_j
    \hat{f}
    \Big{]}
    }_p
    & 
    =
    e^{-t} 
    \bignorm{
    \F^{-1} \Big{[}
    e^{ - \frac{4t}{|\x|^2} \bigc{ 1 \mp \sqrt{1 - \frac{4}{|\x|^2}} }^{-2} }
    \hp_j
    \hat{f}
    \Big{]}
    }_p
    \\
    & \leq e^{-t} \label{3.94}
    \bignorm{
    \F^{-1} \Big{[}
    e^{ - \frac{4t}{|\x|^2} \bigc{ 1 \mp \sqrt{1 - \frac{4}{|\x|^2}} }^{-2} }
    \hp_j
    \Big{]}
    }_1
    \norm{  f}_p.
    \end{align}
    We note that, for all $j \geq 3$,
    \begin{align*}
        \bignorm{
        \F^{-1} \Big{[}
        e^{ - \frac{4t}{|\x|^2} \bigc{ 1 + \sqrt{1 - \frac{4}{|\x|^2}} }^{-2} }
        \hp_j
        \Big{]}
        }_1
        & = 
        \bignorm{
        \F^{-1} \Big{[}
        e^{ - \frac{4t}{ 2^{2j} |\x|^2} \bigc{ 1 + \sqrt{1 - \frac{4}{ 2^{2j} |\x|^2}} }^{-2} }
        \hp_0
        \Big{]}
        }_1
        \\
        & \leq C
        \bignorm{
        e^{ - \frac{4t}{ 2^{2j} |\x|^2} \bigc{ 1 + \sqrt{1 - \frac{4}{ 2^{2j} |\x|^2}} }^{-2} }
        \hp_0
        }_{W^{2,2}}
        \\
        & \leq C. 
        &\hfill \qed
    \end{align*}

We now discuss the optimality of the estimates in the previous section. In particular, we will prove that the low-frequency bound obtained is optimal.
First, we 
denote 
\begin{align*}
    \x_{t} & \coloneqq (\x_1, t^{+1/4}\x_2, t^{+1/4}\x_3), \\
    \x_{t^{-1}} & \coloneqq (\x_1, t^{-1/4}\x_2, t^{-1/4}\x_3).
\end{align*}
We also take a nonnegative nonzero function $\hps \in C^\infty_0$ such that
\[
\supp \hps \subseteq \{ \x\in\R^3 \ | \ |\x|\in({1}/{2} , 1), \  |\x_1| \geq 1/2 \}, \quad
\hps(-\x) = \hps(\x), \text{ for all } \x \in \R^3.
\]

\begin{prop}(Optimality of Low-Frequency Linear Estimate) \label{bnd from below}
There exists a constant $C$ such that, for all $t$ sufficiently large, 
\begin{align*}
    \bignorm{ \F^{-1} \Bigf{ e^{t\lambda_\pm}
    \hps(t^{1/2}\x_{t}) 
    } }_\infty \geq Ct^{-2}.
\end{align*}
\end{prop}

\begin{rem}

We note that this bound from below on the low-frequency estimate is sufficient to prove that for all $t$ sufficiently large, 
\begin{align*}
    \bignorm{ \F^{-1} \Bigf{ \sum_{j\leq 2} \hp_j e^{t \lambda_\pm} } }_\infty \geq Ct^{-2}.
\end{align*}
Indeed, for all $t\geq 1$, we get by a simple application of Young's convolution inequality:
\begin{align*}
    \bignorm{ \F^{-1} \Bigf{ e^{t \lambda_\pm} 
    \hps(t^{1/2}\x_{t})
    } }_\infty
    & \leq
    \bignorm{ \F^{-1} \Bigf{ \hps(t^{1/2}\x_{t}) } }_1
    \bignorm{ \F^{-1} \Bigf{ \sum_{j\leq 2} \hp_j e^{t \lambda_\pm} } }_\infty
    \\
    & =
    \bignorm{ \F^{-1} \Bigf{ \hps(\x) } }_1
    \bignorm{ \F^{-1} \Bigf{ \sum_{j\leq 2} \hp_j e^{t \lambda_\pm} } }_\infty
    \\
    & \leq
    \norm{  \hps(\x) }_{W^{2,2}}
    \bignorm{ \F^{-1} \Bigf{ \sum_{j\leq 2} \hp_j e^{t \lambda_\pm} } }_\infty
    \\
    & \leq 
    C \bignorm{ \F^{-1} \Bigf{ \sum_{j\leq 2} \hp_j e^{t \lambda_\pm} } }_\infty.
\end{align*}
\end{rem}

\begin{pf}
First, note that
\begin{align*}
    &
    \sup_{x\in\R^3} \Big{|} 
    \intr
    e^{ix\cdot\x}
    e^{\pm it \frac{|\x|^2}{2} \sqrt{ \frac{4}{|\x|^2} - 1 }}
    e^{-t\frac{|\x|^2}{2}}
    \hps ( t^{1/2}\x_{t} )
    \dd\x
    \Big{|}
    \\
    &
    \geq t^{-3/2} 
    \Big{|} \intr
    e^{ \pm  i t^{1/2}  \bigc{ \x_1 + |\x|\sqrt{ 1 - \frac{|\x|^2}{4t} }}}
    e^{-\frac{|\x|^2}{2}}
    \hps ( \x_{t} )
    \dd\x
    \Big{|}
    \\
    &
    \geq
    t^{-3/2} 
    \Big{|} \int_{\R^3 \cap \{ \x_1 < 0 \} }
    e^{ \pm  i t^{1/2}  \bigc{ \x_1 + |\x|\sqrt{ 1 - \frac{|\x|^2}{4t} }}}
    e^{-\frac{|\x|^2}{2}}
    \hps ( \x_{t} )
    \dd\x
    \Big{|}
    \\
    &
    \quad -
    t^{-3/2} 
    \Big{|} \int_{ \R^3 \cap \{ \x_1 > 0 \} }
    e^{ \pm  i t^{1/2}  \bigc{ \x_1 + |\x|\sqrt{ 1 - \frac{|\x|^2}{4t} }}}
    e^{-\frac{|\x|^2}{2}}
    \hps ( \x_{t} )
    \dd\x
    \Big{|}
    ,
\end{align*}
where in the above, we chose the specific value $x = (\pm t,0,0)$, and have split the integral into two parts: one with $\x_1 < 0$, which we shall bound from below; and one with $\x_1 > 0$, which will decay at a slightly faster rate.

We consider the $\x_1 < 0$ part. The key step is rewriting the exponent of our oscillation function. Note that in this case, $\x_1 = -|\x_1|$.
\begin{align*}
    i t^{1/2}  \bigc{ \x_1 + |\x|\sqrt{ 1 - {|\x|^2}/{4t} }} 
    & =
    i t^{1/2}  \Bigc{ \x_1 + |\x| - \frac{|\x|^3}{ 4t + 4t \sqrt{ 1-{|\x|^2}/{4t} } } }
    \\
    & =
    i t^{1/2}  \Bigc{ \frac{ -\x_1^2 + |\x|^2 }{|\x_1| + |\x|} - \frac{|\x|^3}{ 4t + 4t \sqrt{ 1-{|\x|^2}/{4t} } } }.
\end{align*}
Then, by the substitution $\x \to \x_{t^{-1}}$, we obtain
\begin{align*}
    & \Big{|} \int_{ \R^3 \cap \{ \x_1 < 0 \} }
    e^{ \pm  i t^{1/2}  \bigc{ \x_1 + |\x|\sqrt{ 1 - \frac{|\x|^2}{4t} }}}
    e^{-\frac{|\x|^2}{2}}
    \hps ( \x_{t} )
    \dd\x
    \Big{|}
    \\
    & = t^{-1/2}
    \Big{|} \int_{ \R^3 \cap \{ \x_1 < 0 \} }
    e^{ \pm  i t^{1/2}  \Bigc{ \frac{ -\x_1^2 + |\x_{t^{-1}}|^2 }{|\x_1| + |\x_{t^{-1}}|} - \frac{|\x_{t^{-1}}|^3}{ 4t + 4t \sqrt{ 1-{|\x_{t^{-1}}|^2}/{4t} } } }}
    e^{-\frac{|\x_{t^{-1}}|^2}{2}}
    \hps ( \x )
    \dd\x
    \Big{|},
\end{align*}
and thus, by the Lebesgue dominated convergence theorem, we have the following convergence as $t\to\infty$:
\begin{align*}
    & t^{+1/2}
    \Big{|} \int_{ \R^3 \cap \{ \x_1 < 0 \} }
    e^{ \pm  i t^{1/2}  \bigc{ \x_1 + |\x|\sqrt{ 1 - \frac{|\x|^2}{4t} }}}
    e^{-\frac{|\x|^2}{2}}
    \hps ( \x_{t} )
    \dd\x
    \Big{|}
    \\
    & \to
    \Big{|} \int_{ \R^3 \cap \{ \x_1 < 0 \} }
    e^{ \pm  i { (\x_2^2 + \x_3^2) }/{2|\x_1|} }
    e^{-\frac{\x_1^2}{2}}
    \hps ( \x )
    \dd\x
    \Big{|}, \text{ as } t\to\infty,
\end{align*}
where the final integral is a positive constant, by the restriction to the support of $\hps$.
We thus conclude that there exists some constant C such that, for all sufficiently large $t$,
\begin{align*}
    \Big{|} \int_{ \R^3 \cap \{ \x_1 < 0 \} }
    e^{ \pm  i t^{1/2}  \bigc{ \x_1 + |\x|\sqrt{ 1 - \frac{|\x|^2}{4t} }}}
    e^{-\frac{|\x|^2}{2}}
    \hps ( \x_{t} )
    \dd\x
    \Big{|}
    \geq C t^{-1/2}.
\end{align*}

It remains to show that the $\x_1 >0$ part decays at a faster rate. The proof is similar to that of Proposition \ref{appendix prop 1} in the Appendix, and so we skip the details. 
The key difference is that $\x_1>0$ prevents the cancellation of $\x_1$ in the exponent which we exploited above for the $\x_1 <0$ case. In fact, in the $\x_1>0$ case, we may extract an extra decay of $t^{-1/4}$ compared to the normal bound from above proven in Proposition \ref{appendix prop 1}. (See also Remark~\ref{last remark}).
\end{pf}
\bigskip

We lastly must find some suitable initial data $(u_0, a_0)$ such that 
\begin{align*}
    \bignorm{ \F^{-1} 
    \begin{bmatrix}
        \ha(t) \\
        \hv(t)
    \end{bmatrix} 
    }_\infty \geq C t^{-2}.
\end{align*}
Recall that
\[
v_0 = |D|^{-1} \div ( u_0 ).
\]

\begin{prop} \label{proof of 3rd}
Let 
\[
u_0 \coloneqq e^{-|x|^2} \begin{bmatrix}
    1 \\
    1 \\
    1
\end{bmatrix},
\]
and let $a_0 \in L^1 \cap \dot{B}^{0}_{\infty,1}$.
Then there exist constants $C,\delta>0$ such that, if $\norm{a_0}_{L^1 \cap \dot{B}^{0}_{\infty,1}} \leq \delta$, then for all sufficiently large $t$,
\begin{align*}
    \bignorm{\F^{-1}
    \begin{bmatrix} 
        \ha(t) \\
        \hv(t)
    \end{bmatrix} 
    }_\infty \geq C t^{-2}.
\end{align*}
\end{prop}

\begin{pf}
Firstly, we note that the smallness of $a_0$ is required to distinguish between $v_0$ and $a_0$ via the following reverse triangle inequality
\begin{align*}
    \bignorm{ \F^{-1} 
    \begin{bmatrix}
        \ha(t) \\
        \hv(t)
    \end{bmatrix} 
    }_\infty 
    & \geq
    \bignorm{ \F^{-1} 
    \Big{[}
        \ha(t) 
    \Big{]}
    }_\infty 
    \\
    & = 
    \bignorm{ \F^{-1} 
    \Big{[}
    \frac{ e^{t\lambda_-} \lambda_+ - e^{t\lambda_+} \lambda_- }{ \lambda_+ - \lambda_- }    
    \ha_0 
    + 
    \frac{ \bigc{ e^{t\lambda_-} - e^{t\lambda_+}  } \lambda_+ }{ \lambda_+ - \lambda_- }
    \hv_0
    \Big{]}
    }_\infty 
    \\
    & \geq
    \bignorm{ \F^{-1} 
    \Big{[}
        \frac{ \bigc{ e^{t\lambda_-} - e^{t\lambda_+}  }\lambda_+ }{ \lambda_+ - \lambda_- }
    \hv_0
    \Big{]}
    }_\infty 
    -
    C \delta t^{-2}.
\end{align*}
Next, focusing on the norm with $\hat{v}_0$, we note that, due to our choice of $u_0$,
\begin{align*}
    \hv_0 = \frac{i C_1 ( \x_1 + \x_2 + \x_3 ) e^{-C_2 |\x|^2} }{|\x|}.
\end{align*}
Now, we inspect the $L^\infty$-norm above. As in the proof of the previous proposition, we insert the function $\hps$ to restrict our function to the low-frequency case.
\begin{align*}
    &  \bignorm{ \F^{-1} 
    \Big{[}
        \frac{ \bigc{ e^{t\lambda_-} - e^{t\lambda_+}  } \lambda_+ }{ \lambda_+ - \lambda_- }
    \hv_0
    \Big{]}
    }_\infty
    \geq C
    \bignorm{ \F^{-1} 
    \Big{[}
        \frac{ \bigc{ e^{t\lambda_-} - e^{t\lambda_+}  } \lambda_+ }{ \lambda_+ - \lambda_- }
    \hv_0
        \hps (t^{1/2}\x_t)
    \Big{]}
    }_\infty
    \\
    &  = C
    \sup_{x\in \R^3} \Big{|} 
    \intr 
    e^{ix\cdot\x} 
    \hv_0(\x) 
    \hps (t^{1/2}\x_t)
    \frac{ 
    \Bigc{
    e^{-t\frac{|\x|^2}{2} \bigc{1 - i \sqrt{ \frac{4}{|\x|^2} - 1 }} } 
    -  
    e^{-t\frac{|\x|^2}{2} \bigc{1 + i \sqrt{ \frac{4}{|\x|^2} - 1 }} } 
    }
    }{2-(|\x|^2/2)+(i|\x|/2)\sqrt{ 4-|\x|^2 }} 
    \dd\x
    \Big{|}
    \\
    & \geq
    C
    t^{-3/2} \Big{|} 
    \intr 
    e^{-\frac{|\x|^2}{2}} 
    \hv_0(t^{-1/2}\x) 
    \hps (\x_t)
    \frac{ 
    \Bigc{
    e^{-it^{1/2} \bigc{\x_1 - |\x| \sqrt{ 1 - \frac{|\x|^2}{4t} }} } 
    -  
    e^{-it^{1/2} \bigc{\x_1 + |\x| \sqrt{ 1 - \frac{|\x|^2}{4t} }} } 
    }
    }{2-(|\x|^2/2t)+(i|\x|/2t^{1/2})\sqrt{ 4-t^{-1}|\x|^2 }} 
    \dd\x
    \Big{|}.
\end{align*}
Where, for the last step, we set $x=(-t,0,0)$ and made the substitution $\x\to t^{-1/2}\x$, as we did in the previous proposition. 
The differing signs in our two exponential functions must be addressed. We separate the two terms into two integrals, and make the substitution $\x \to -\x$ for one of them. Then, since our $\hv_0$ is antisymmetric, we get
\begin{align*}
    & \Big{|} 
    \intr 
    e^{- \frac{|\x|^2}{2}} 
    \hv_0(t^{-1/2}\x) 
    \hps ( \x_t )
    \frac{ 
    \Bigc{
    e^{-it^{1/2} \bigc{\x_1 - |\x| \sqrt{ 1 - \frac{|\x|^2}{4t} }} } 
    -  
    e^{-it^{1/2} \bigc{\x_1 + |\x| \sqrt{ 1 - \frac{|\x|^2}{4t} }} } 
    }
    }{2-(|\x|^2/2t)+(i|\x|/2t^{1/2})\sqrt{ 4-t^{-1}|\x|^2 }} 
    \dd\x
    \Big{|}
    \\
    & = 
    \Big{|} 
    \intr 
    e^{- \frac{|\x|^2}{2}} 
    \hv_0(t^{-1/2}\x) 
    { \hps ( \x_t )}
    \frac{ 
    \Bigc{
    e^{+it^{1/2} \bigc{\x_1 + |\x| \sqrt{ 1 - \frac{|\x|^2}{4t} }} } 
    +
    e^{-it^{1/2} \bigc{\x_1 + |\x| \sqrt{ 1 - \frac{|\x|^2}{4t} }} } 
    }
    }{2-(|\x|^2/2t)+(i|\x|/2t^{1/2})\sqrt{ 4-t^{-1}|\x|^2 }} 
    \dd\x
    \Big{|}
    \\
    & = 
    2 \Big{|} 
    \intr 
    e^{- \frac{|\x|^2}{2}} 
    \hv_0(t^{-1/2}\x) 
    { \hps ( \x_t )}
    \frac{ 
    \cos( t^{1/2} \bigc{\x_1 + |\x| \sqrt{ 1 - {|\x|^2}/{4t} }} )
    }{2-(|\x|^2/2t)+(i|\x|/2t^{1/2})\sqrt{ 4-t^{-1}|\x|^2 }} 
    \dd\x
    \Big{|}.
\end{align*}
We may then follow the same steps as the previous proof, as $\cos( \cdot )$ is sufficiently similar to $e^{i \cdot }$ for our purposes. 
\end{pf}
\bigskip

\noindent\textbf{Proof of Theorem 1.2.} 
\eqref{thm1 l1 linf} is proven by Proposition \ref{proof of 1st} and Proposition \ref{prop no singularity}. 
Indeed, we have for all $s \in \R$, $p \in [2,\infty],$ and $q \in [1,\infty]$,
\begin{align*}
        & \bignorm{
        \begin{bmatrix}
        a(t) \\
        v(t)
        \end{bmatrix}  }_{\dot{B}^s_{p,q}}
        \\
        &
        =
        \Biggc{
        \sum_{j\leq-1} 2^{sjq} \bignorm{ 
        \ddj 
        \begin{bmatrix}
        a(t) \\
        v(t)
        \end{bmatrix} 
        }_p^q
        +
        \sum_{j = 0}^{2} 2^{sjq} \bignorm{ 
        \ddj 
        \begin{bmatrix}
        a(t) \\
        v(t)
        \end{bmatrix} 
        }_p^q
        +
        \sum_{j\geq3} 2^{sjq} \bignorm{ 
        \ddj 
        \begin{bmatrix}
        a(t) \\
        v(t)
        \end{bmatrix} 
        }_p^q
        }^{1/q}
        \\
        &
        \leq 
        C t^{-\frac{3}{2} (1-\frac{1}{p}) - \frac{1}{2} (1 - \frac{2}{p}) }
        \bignorm{ \F^{-1} 
        \begin{bmatrix}
        \ha_0 \\
        \hv_0
        \end{bmatrix}  }^l_{\dot{B}^{s}_{1,q}}
        +
        C e^{-t}
        \bignorm{ \F^{-1} 
        \begin{bmatrix}
        \ha_0 \\
        \hv_0
        \end{bmatrix}  }^h_{
        \dot{B}^{s+3(1-1/p)}_{1,q} 
        }.
\end{align*}

\eqref{thm1 lp lp} is proven by Proposition \ref{proof of 2nd}. 
Indeed, there exists a constant $C$ such that, for all $p \in [1,\infty]$ and $t>0$,
\begin{align*}
    \bignorm{
    \F^{-1} \Big{[}
    e^{t\lambda_\pm}
    \hp_j
        \begin{bmatrix} 
        \ha_0 \\
        \hv_0
    \end{bmatrix} 
    \Big{]}
    }_p
    & 
    \leq
    C e^{-t} 
    \bignorm{
    \F^{-1} \Big{[}
    (\hp_{j-1} + \hp_j + \hp_{j+1})
        \begin{bmatrix} 
        \ha_0 \\
        \hv_0
    \end{bmatrix} 
    \Big{]}
    }_p
    \text{ for all } 
    j \geq  3.
\end{align*}
We thus obtain
\begin{align*}
        \bignorm{ \F^{-1} 
        \begin{bmatrix}
        \ha(t) \\
        \hv(t)
        \end{bmatrix}  }^h_{\dot{B}^s_{p,q}}
        \leq 
        C e^{-t}
        \bignorm{ \F^{-1} 
        \begin{bmatrix}
        \ha_0 \\
        \hv_0
        \end{bmatrix}  }^h_{\dot{B}^{s}_{p,q} }.
\end{align*}

\eqref{thm 1 linf below} is proven by Proposition \ref{proof of 3rd}. \hfill\qed

\appendix
\section{Behaviour of Solutions near $|\x| =2$}\label{appendix no singularity}

\begin{prop}\label{prop no singularity}
    Let $p \in \{ 2, \infty \}$, $j \in \{ 0,1,2 \}$. Then there exist constants $C,c>0$ such that, for all $t>0$,
    \begin{align*}
    \bignorm{  
    \ddj
        \begin{bmatrix}
        a(t) \\
        v(t)
        \end{bmatrix}
    }_{p}
    \leq
    Ce^{-ct}.
    \end{align*}
\end{prop}

\begin{pf}
Taking directly the $L^2$ norm of $\ddj a, \ddj v$, for $j \in \{ 0, 1, 2 \}$, we use Plancherel's theorem to take the Fourier transform inside the norm. Looking at 
\[
\frac{e^{t\lambda_-}\lambda_+ - e^{t\lambda_+ } \lambda_- }{ \lambda_+ - \lambda_- } 
\]
and the other similar expressions in \eqref{a formula} and \eqref{v formula}, we see that the numerator and denominator both approach $0$ as $|\x|\to 2$. 
First, we may rewrite
\[
\frac{e^{t\lambda_-}\lambda_+ - e^{t\lambda_+ } \lambda_- }{ \lambda_+ - \lambda_- } 
=
\frac{(e^{t\lambda_-} - e^{t\lambda_+ }) \lambda_+ }{ \lambda_+ - \lambda_- }
+
\frac{e^{t\lambda_+} (\lambda_+ - \lambda_-)  }{ \lambda_+ - \lambda_- }.
\]
The second term is of course just $e^{t\lambda_+}$. Focusing on the first term, we inspect the $L^2$ norm around $|\x|=2$.
\begin{align*}
    & 
    \bignorm{ \frac{(e^{t\lambda_-} - e^{t\lambda_+ }) \lambda_+ }{ \lambda_+ - \lambda_- } 
    (\hp_{j-1} + \hp_{j} +\hp_{j+1} )
    }^2_2
    \\
    & 
    \leq
    C
    \int_{ \{ \frac{1}{4} \leq |\x| \leq 2 \} }
    \Bigg{|}
    \frac{ e^{- t \frac{|\x|^2}{2} } 
    \bigc{
    e^{+it\frac{|\x|^2}{2} \sqrt{ \frac{4}{|\x|^2} - 1 } } 
    -
    e^{-it\frac{|\x|^2}{2} \sqrt{ \frac{4}{|\x|^2} - 1 } } 
    }
    }{\sqrt{ 4/|\x|^2 - 1 }  }
    \Bigg{|}^2
    \dd\x
    \\
    & \quad \quad \quad \quad \quad \quad \quad \quad 
    + C
    \int_{ \{ 2 \leq |\x| \leq 8 \} }
    \Bigg{|}
    \frac{
    e^{- t \frac{|\x|^2}{2} } 
    \bigc{
    e^{+t\frac{|\x|^2}{2} \sqrt{ 1 - \frac{4}{|\x|^2} } } 
    -
    e^{-t\frac{|\x|^2}{2} \sqrt{ 1 - \frac{4}{|\x|^2} } }} 
    }
    {\sqrt{ 1 -  4/|\x|^2  }  }
    \Bigg{|}^2
    \dd\x
    \\
    & 
    \eqqcolon I_1 + I_2.
\end{align*}
Looking at $I_1$ first, for all $t>0$, we have
\begin{align*}
    I_1 & \leq C e^{-t/32}
    \int_{ \{ \frac{1}{4} \leq |\x| \leq 2 \} }
    \Bigg{|}
    \frac{  
    \bigc{
    e^{+it\frac{|\x|^2}{2} \sqrt{ \frac{4}{|\x|^2} - 1 } } 
    -
    e^{-it\frac{|\x|^2}{2} \sqrt{ \frac{4}{|\x|^2} - 1 } } 
    }
    }{\sqrt{ 4/|\x|^2 - 1 }  }
    \Bigg{|}^2
    \dd\x
    \\
    & \leq C t^2 e^{-t/32} 
    \int_{ \{ \frac{1}{4} \leq |\x| \leq 2 \} }
    \Bigg{|}
    \frac{  
    \sin(t\frac{|\x|^2}{2} \sqrt{ \frac{4}{|\x|^2} - 1 } )
    }{ t\frac{|\x|^2}{2}  \sqrt{ 4/|\x|^2 - 1 }  }
    \Bigg{|}^2
    \dd\x
    \\
    & \leq C e^{-ct}.
\end{align*}
For $I_2$, values of $|\x|$ away from 2 are handled as in Proposition \ref{mainprop}. It remains only to consider the integral over $|\x| \in (2,2+\e)$, for some small $\e>0$.
To that end, we first extract our time decay
\begin{align*}
    & \int_{ \{ 2 \leq |\x| \leq 2+\e \} }
    \Bigg{|}
    \frac{
    e^{- t \frac{|\x|^2}{2} } 
    \bigc{
    e^{+t\frac{|\x|^2}{2} \sqrt{ 1 - \frac{4}{|\x|^2} } } 
    -
    e^{-t\frac{|\x|^2}{2} \sqrt{ 1 - \frac{4}{|\x|^2} } }
    } 
    }
    {\sqrt{ 1 -  4/|\x|^2  }  }
    \Bigg{|}^2
    \dd\x
    \\
    & \leq 
    e^{ -4t }
    \int_{ \{ 2 \leq |\x| \leq 2+\e \} }
    \Bigg{|}
    \frac{
    e^{+t\frac{|\x|^2}{2} \sqrt{ 1 - \frac{4}{|\x|^2} } } 
    -
    e^{-t\frac{|\x|^2}{2} \sqrt{ 1 - \frac{4}{|\x|^2} } } 
    }
    {\sqrt{ 1 -  4/|\x|^2  }  }
    \Bigg{|}^2
    \dd\x,
\end{align*}
and consider the Taylor expansion of the remaining integrand:
\begin{align*}
    & \frac{
    e^{+t\frac{|\x|^2}{2} \sqrt{ 1 - \frac{4}{|\x|^2} } } 
    -
    e^{-t\frac{|\x|^2}{2} \sqrt{ 1 - \frac{4}{|\x|^2} } }
    }{\sqrt{ 1 - 4/|\x|^2 }  }
    \\
    &
    =
    \frac{1}{\sqrt{ 1 - 4/|\x|^2 }  } \sum^{\infty}_{n=1}
    \frac{ (t\frac{|\x|^2}{2})^n \bigc{  (1 - 4/|\x|^2)^{n/2} - (-1)^n (1 - 4/|\x|^2)^{n/2} } }{n!}
    \\
    & =
    \frac{ 2 }{\sqrt{ 1 - 4/|\x|^2 }  } \sum_{n \geq 1, \text{odd} }
    \frac{ (t\frac{|\x|^2}{2})^n (1 - 4/|\x|^2)^{n/2} }{n!}
    \\
    & \leq t|\x|^2 \sum_{n \geq 0, \text{even} }
    \frac{ (t\frac{|\x|^2}{2})^n (1 - 4/|\x|^2)^{n/2} }{n!}
    \\
    & \leq t |\x|^2 e^{    +t\frac{|\x|^2}{2} \sqrt{ 1 - \frac{4}{|\x|^2} }     } ,
\end{align*}
for all $|\x| > 2$, $t>0$.
We thus obtain for all $t>0$ and any $\e$ sufficiently small
\begin{align*}
    I_2 & 
    \leq C e^{-4t}
    \int_{ \{ 2 \leq |\x| \leq 2+\e \} }
    \Big{|} 
     t |\x|^2 e^{    +t\frac{|\x|^2}{2} \sqrt{ 1 - \frac{4}{|\x|^2} }     } 
    \Big{|}^2
    \dd\x
    + 
    C e^{-ct}
    \\
    & \leq C e^{-ct},
\end{align*}
and so we have obtained the result for $p=2$.

Using the inequality
\[
\norm{ \F^{-1} [\hat{f}] }_\infty
\leq
\norm{ \hat{f} }_1,
\]
the same result is proven similarly for the $L^\infty$ norm.
\end{pf}


\section{Bounds from Above for Low Frequencies}\label{appendix prop 1}

\begin{prop}
Let $j\leq -1$. Then, for all sufficiently large $t$,
    \begin{align*}
    \sup_{x\in\R^3} \Big{|} 
    \intr 
    e^{i (t^{-1/2}x)\cdot\x}
    e^{\pm i t^{1/2} \frac{|\x|^2}{2} \sqrt{ \frac{4}{|\x|^2} - t^{-1} }}
    e^{-\frac{|\x|^2}{2}}
    \hp(2^{-j} t^{-1/2} \x )
    \dd\x
    \Big{|} 
    & \leq C t^{- 1/2}.
\end{align*}
\end{prop}

\begin{pf}
First note that, by radial symmetry, we can assume without loss of generality that $x = (x_1,0,0)$.
We then split the frequency space into four subsets:
\begin{align*}
    B_1 & \coloneqq
    \{
    \x \in \R^3 \ | \ |\x_2| , |\x_3| \leq t^{-1/4}
    \} \\
    B_2 & \coloneqq
    \{
    \x \in \R^3 \ | \ |\x_2| \geq t^{-1/4} , |\x_3| \leq t^{-1/4}
    \}
    \\
    B_3 & \coloneqq
    \{
    \x \in \R^3 \ | \ |\x_2| \leq t^{-1/4}, |\x_3| \geq t^{-1/4}
    \}
    \\
    B_4 & \coloneqq
    \{
    \x \in \R^3 \ | \  |\x_2| , |\x_3| \geq t^{-1/4}
    \} 
    \\
    \R^3 & =  
    B_1 \cup B_2 \cup B_3 \cup B_4.
\end{align*}
We shall estimate our integral over each of the four subsets. We start with $B_1$.
\begin{align*}
    &
    \sup_{x\in\R^3}
    \Big{|} \int_{\x \in B_1}
    e^{i t^{-1/2} x_1 \x_1}
    e^{\pm i t^{1/2} \frac{|\x|^2}{2} \sqrt{ \frac{4}{|\x|^2} - t^{-1} }}
    e^{-\frac{|\x|^2}{2}}
    \hp(2^{-j} t^{-1/2} \x )
    \dd\x
    \Big{|}
    \\
    & \leq 
    \int_{ \R }
    \int_{ \{ |\x_3| \leq t^{-1/4} \} } 
    \int_{ \{ |\x_2| \leq t^{-1/4} \} }
    e^{-\frac{|\x|^2}{2}}
    \dd\x_2 \dd\x_3 \dd\x_1
    \leq C t^{-1/2}.
\end{align*}

For $B_2,$ we first note that
\begin{align*}
    e^{ i t^{1/2} \frac{|\x|^2}{2} \sqrt{ \frac{4}{|\x|^2} - t^{-1} } }
    =
    \Bigc{ \pt_j e^{ i t^{1/2} \frac{|\x|^2}{2} \sqrt{ \frac{4}{|\x|^2} - t^{-1} } } }
    \frac{ \frac{|\x|^2}{2} \sqrt{ \frac{4}{|\x|^2} - t^{-1} } }{ 2 \x_j i t^{1/2} (2 - |\x|^2 t^{-1}) },
\end{align*}
and thus,
\begin{align*}
    & 
    \sup_{x\in\R^3}
    \Big{|} \int_{\x \in B_2}
    e^{i t^{-1/2} x_1 \x_1}
    e^{\pm i t^{1/2} \frac{|\x|^2}{2} \sqrt{ \frac{4}{|\x|^2} - t^{-1} }}
    e^{-\frac{|\x|^2}{2}}
    \hp(2^{-j} t^{-1/2} \x )
    \dd\x
    \Big{|}
    \\
    &  \leq 
    \int_{ \{ |\x_1| < t^{1/2} \} } \int_{ \{ |\x_3| \leq t^{-1/4} \} } 
    \\
    &
    \quad \quad \quad \Big{|} 
    \int_{ \{ |\x_2| \geq t^{-1/4} \} }
    \frac{ 
    \Bigc{ 
    \pt_2 e^{ i t^{1/2} \frac{|\x|^2}{2} \sqrt{ \frac{4}{|\x|^2} - t^{-1} } } 
    }
    \frac{|\x|^2}{2} \sqrt{ \frac{4}{|\x|^2} - t^{-1} } 
    \hp(2^{-j} t^{-1/2} \x )
    e^{ - \frac{|\x|^2}{2}}
    }
    { 
    2 \x_2 i t^{1/2} (2 - |\x|^2 t^{-1}) 
    } 
    \dd\x_2
    \Big{|}
    \dd\x_3
    \dd\x_1
    \\
    & \leq 
    \int_{ \{ |\x_1| < t^{1/2} \} } \int_{ \{ |\x_3| \leq t^{-1/4} \} }
    \\
    & 
    \quad \quad \quad 
    \Biggc{
    Ct^{-1/4} + Ct^{-1/2}
    \int^{\infty}_{t^{-1/4}} 
    \Big{|}
    \pt_2 \Bigc{ 
    \frac{ 
    \frac{|\x|^2}{2} \sqrt{ \frac{4}{|\x|^2} - t^{-1} } 
    \hp(2^{-j} t^{-1/2} \x ) 
    e^{ - \frac{|\x|^2}{2}}
    }
    { 2 \x_2 i  (2 - |\x|^2 t^{-1}) }
    }
    \Big{|}
    \dd\x_2
    }
    \dd\x_3
    \dd\x_1.
\end{align*}
The above is accomplished by a simple integration by parts. Routine checks verify that
\begin{align*}
    \int^{\infty}_{t^{-1/4}}
    \Big{|}
    \pt_2 \Bigc{ \frac{ \frac{|\x|^2}{2} \sqrt{ \frac{4}{|\x|^2} - t^{-1} } }{ 2 \x_2 i  (2 - |\x|^2 t^{-1}) } 
    \hp(2^{-j} t^{-1/2} \x ) 
    e^{ - \frac{\x_2^2}{2}}
    }
    \Big{|}
    \dd\x_2
    \leq Ct^{+1/4}.
\end{align*}
An additional time decay of $t^{-1/4}$ is obtained from the integral over $\x_3$, and boundedness of the integral over $\x_1$ is assured by the presence of the Gauss function $e^{-\x_1^2/2}$. Thus 
\begin{align*}
    \sup_{x\in\R^3}
    \Big{|} \int_{\x \in B_2}
    e^{i t^{-1/2} x_1 \x_1}
    e^{\pm i t^{1/2} \frac{|\x|^2}{2} \sqrt{ \frac{4}{|\x|^2} - t^{-1} }}
    e^{-\frac{|\x|^2}{2}}
    \hp(2^{-j} t^{-1/2} \x )
    \dd\x
    \Big{|}
    \leq C t^{-1/2}.
\end{align*}

The process for bounding the integral over $B_3$ is identical to that of $B_2$. We are thus left with the integral over $B_4$.

We note that 
\begin{align*}
    & e^{ i t^{1/2} \frac{|\x|^2}{2} \sqrt{ \frac{4}{|\x|^2} - t^{-1} } }
    \\
    & =
    \frac{\pt_3 \pt_2 e^{ i t^{1/2} \frac{|\x|^2}{2} \sqrt{ \frac{4}{|\x|^2} - t^{-1} } }}{4\x_2 \x_3 i t^{1/2}} 
    \frac{\frac{|\x|^6}{2} \sqrt{  \frac{4}{|\x|^2} - t^{-1}  }^3 }{ 4 i t^{1/2} |\x|^2  \sqrt{ \frac{4}{|\x|^2} - t^{-1} } - it^{-1/2} |\x|^6 \sqrt{  \frac{4}{|\x|^2} - t^{-1}  }^3 - 4 },
    \\
    & \eqqcolon (A)
\end{align*}
and thus, by integration by parts,
\begin{align*}
    &
    \sup_{x\in\R^3}
    \Big{|} \int_{\x \in B_4}
    e^{i t^{-1/2} x_1 \x_1}
    e^{\pm i t^{1/2} \frac{|\x|^2}{2} \sqrt{ \frac{4}{|\x|^2} - t^{-1} }}
    e^{-\frac{|\x|^2}{2}}
    \hp(2^{-j} t^{-1/2} \x )
    \dd\x
    \Big{|}
    \\
    & \leq 
    \int^{1}_{-1}
    \Big{|}
    \int_{ \{ |\x_3| \geq t^{-1/4} \} } 
    \int_{ \{ |\x_2| \geq t^{-1/4} \} }
    (A) \, \hp(2^{-j} t^{-1/2} \x ) e^{-\frac{|\x|^2}{2}}
    \dd\x_2
    \dd\x_3
    \Big{|}
    \dd\x_1
    \\
    & \leq 
    Ct^{-1/2},
\end{align*}
by some longer but still routine calculations similar to the integral over $B_2$.
Therefore, we have obtained for all $j \leq -1,$ and all $t > t_0$,
\begin{align*}
    & \sup_{x\in\R^3} 
    \Big{|}
    \intr
    e^{ix\cdot\x}
    e^{\pm it\frac{|\x|^2}{2} \sqrt{ \frac{4}{|\x|^2} - 1 }}
    e^{-t\frac{|\x|^2}{2}}
    \hat{\phi}(2^{-j} \x)
    \dd\x
    \Big{|}
    \leq C t^{-2}.
\end{align*}
This completes the proof.
\end{pf}

\begin{rem}\label{last remark}
    When we limit ourselves to the $\x_1 > 0$ case, as in the proof of Proposition \ref{bnd from below}, we may obtain an extra decay of $t^{-1/4}$ by performing an additional integration by parts over $\x_1$. 
\end{rem}







\bibliographystyle{unsrt} \bibliography{sample.bib}

\begin{bibdiv}
\begin{biblist}

\bib{Danchinbook}{book}{
      author={Bahouri, Hajer},
      author={Chemin, Jean-Yves},
      author={Danchin, Rapha\"{e}l},
       title={Fourier analysis and nonlinear partial differential equations},
      series={Grundlehren der mathematischen Wissenschaften [Fundamental
  Principles of Mathematical Sciences]},
   publisher={Springer, Heidelberg},
        date={2011},
      volume={343},
        ISBN={978-3-642-16829-1},
      review={\MR{2768550}},
}

\bib{Danchin2000}{article}{
      author={Danchin, R},
       title={Global existence in critical spaces for compressible
  navier-stokes equations},
        date={2000},
     journal={Invent. math.},
      volume={141},
      number={3},
       pages={579\ndash 614},
}

\bib{Danchin-Xu2017}{article}{
      author={Danchin, R},
      author={Xu, J},
       title={Optimal time-decay estimates for the compressible navier–stokes
  equations in the critical lp framework},
        date={2017},
     journal={Arch. Rational Mech. Anal.},
      volume={224},
      number={3},
       pages={53 \ndash  90},
}

\bib{hoff-zumbrun}{article}{
      author={Hoff, David},
      author={Zumbrun, Kevin},
       title={Multi-dimensional diffusion waves for the {N}avier-{S}tokes
  equations of compressible flow},
        date={1995},
        ISSN={0022-2518},
     journal={Indiana Univ. Math. J.},
      volume={44},
      number={2},
       pages={603\ndash 676},
      review={\MR{1355414}},
}

\bib{kobayashi-shibata}{article}{
      author={Kobayashi, Takayuki},
      author={Shibata, Yoshihiro},
       title={Remark on the rate of decay of solutions to linearized
  compressible {N}avier-{S}tokes equations},
        date={2002},
        ISSN={0030-8730},
     journal={Pacific J. Math.},
      volume={207},
      number={1},
       pages={199\ndash 234},
      review={\MR{1974473}},
}

\bib{matsumura-nishida1979}{article}{
      author={Matsumura, A},
      author={Nishida, T},
       title={The initial value problem for the equation of motion of
  compressible viscous and heat-conductive fluids},
        date={1979},
     journal={Proc. Japan Acad. Ser. A Math. Sci.},
      volume={55},
      number={9},
       pages={337 \ndash  342},
}

\bib{matsumura-nishida1980}{article}{
      author={Matsumura, A},
      author={Nishida, T},
       title={The initial value problem for the equations of motion of viscous
  and heat-conductive gases},
        date={1980},
     journal={J. Math. Kyoto Univ.},
      volume={20},
      number={1},
       pages={67 \ndash  104},
}

\end{biblist}
\end{bibdiv}

\end{document}